\documentclass[12pt]{article}
\newcommand{\End}{\mathop{\rm End}\nolimits}

\newcommand{\Ber}{\mathop{\rm Ber}\nolimits}
\newcommand{\Alt}{\mathop{\rm Alt}\nolimits}
\newcommand{\sgn}{\mathop{\rm sgn}\nolimits}
\newcommand{\diag}{\mathop{\rm diag}\nolimits}

\usepackage{amsfonts}

\title{On the one-dimensional representations 
of the general linear supergroup}

\author{ I.M. Trishin}
\date{}

\begin{document}
\maketitle

Abstract. Because of its multiplicativity,
 the Berezinian  is the character of the 
one-dimensional
representation of the general linear supergroup.
 We give an explicit
 construction of this representation on a
 space of tensors.
Similarly, we construct the representation
such that its character is the inverse
of the Berezinian.

\setcounter{section}{-1}
\section{\rm Introduction}
\indent

 The following result belongs to the  basic 
ones in the classic case.
Let $V$ be an $n$-dimensional linear space over a
 field $K$ of zero characteristic. Then the space
$\Alt (\underbrace{V\otimes\cdots\otimes V}_n)$ 
of the skew-symmetric tensors of degree $n$
is the one-dimensional module 
over the general linear group $GL\, V$.
The character of such representation 
of $GL\, V$ is the 
determinant. In the present paper we construct
a pair of analogous representations
on spaces of tensors for the supercase.
The character of one of these representations 
 is the Berezinian. The other character
is the inverse of the Berezinian.

Outline the way to our construction in general.
Let $G$ be a Grassman algebra with a countable set 
of generators, $V=V_0\oplus V_1$ a free $\mathbb{Z}_2$-graded
$G$-module of dimension $m|n$, where 
$m=\dim_GV_0$, $n=\dim_GV_1$. Suppose $W_{\lambda_h}$,
$W_{\lambda_g}$ are the irreducible $GL\, V$-modules
defined by the partitions 
$\lambda_h=(\underbrace{n+1,\ldots,n+1}_m)$,
$\lambda_g=(\underbrace{n,\ldots,n}_{m+1})$
(respectively)
and $b$ is a formal element such that for 
an arbitrary ${\cal A}\in GL\, V$ we have
${\cal A}(b)=\Ber A\cdot b$,
where $A$ is a matrix of $\cal A$.
We show that the $GL\, V$-modules $W_{\lambda_h}$ 
and $b\cdot W_{\lambda_g}$ are isomorphic
(see Theorem 2.1). 
To prove this theorem we use a natural correspondence
between the sets of the base vectors of $W_{\lambda_h}$
and $W_{\lambda_g}$.
The tensor $\tilde b$ generating
the one-dimensional $GL\, V$-module is constructed
 of  the element $b$ and
these base vectors of $W_{\lambda_h}$ and
$W_{\lambda_g}$ (see Theorem 3.1).
By analogy we construct the 
one-dimensional $GL\, V$-module such that
its character is the inverse of the Berezinian.

This paper is organized  as follows. 
We introduce the main concepts in \S 1.
In \S 2 we prove that 
$W_{\lambda_h}\simeq b\cdot W_{\lambda_g}$
and obtain an important corollary of this result
(see Theorem 2.2). Then
in \S 3 the one-dimensional representations of 
$GL\, V$ are constructed.

In relation to the present paper, 
the result of 
H.M. Khudaverdian and Th.Th. Voronov 
must be mentioned. In  \cite{KV}, 
the Berezinian is expressed as the ratio
of the two  Hankel determinants 
such that these  are the characters of the 
representations of 
$GL\, V$ defined by the partitions $\lambda_h$
and $\lambda_g$.  
On one hand, this is a simple corollary of
Theorem 2.1, on the other  directly implies
this theorem (unfortunately, the last easy 
observation have come too late to be useful for 
the author). 

The author wishes to thank Th.Th. Voronov for useful
discussions.

\section{\rm{The  basic concepts and
  auxiliary results}}
\setcounter{equation}{0}
\indent

 Let $K$ be a  field
 of  characteristic zero,
$G=G_0 \oplus G_1$  the Grassman algebra over 
$K$ with a countable set of generators.

By definition, the notion of $G$-module
includes the property to be free.

Let   $V=V_0\oplus V_1$ be a
finite-dimensional $\mathbb{Z}_2$-graded
 $G$-bimodule and   the structures of 
the left and
 right $G$-modules on $V$  are compatible
  (see, for example, \cite{Det}, \cite{Red}).
Denote
$m=\dim_G V_{0}$, $n=\dim_G V_{1}$.
The pair $m|n$ is the {\it dimension} of $V$.

We suppose that any operation of changing 
base of $V$
 is even.

Let $\End V$ be the algebra of all $G$-linear 
mappings from $V$ to $V$. The algebra $\End V$
 is isomorphic to the algebra
$M=M_0\oplus M_1$ of block matrices with 
Grassman elements. The even component $M_0$ of $M$
is called the {\it  full matrix superalgebra} and 
is denoted
 by $M_{m,n}$. Actually the elements of $M_{m,n}$  
are all matrices of 
the form $A=\left( \begin {array}{l} 
A_{11}\;A_{12}\\A_{21}\;A_{22}
\end{array} \right)$, where $A_{11},\; A_{22}$
are square $G_0$-matrices of orders $m$ and $n$ 
respectively
and $A_{12},\;A_{21}$ are rectangular $G_1$-matrices
of  corresponding orders. 

Let $GL \, V$ and $GL_{m,n}$   be the groups of all 
invertible 
elements of  $\End V$ and $M_{m,n}$ respectively. 
We have the isomorphism
 $GL \, V \simeq GL_{m,n}$. 
The group $GL_{m,n}$ ($GL\, V$)
is 
 the {\it general linear supergroup}. 

Recall that for any
$A=\left( \begin {array}{l} 
A_{11}\;A_{12}\\A_{21}\;A_{22}
\end{array} \right)\in GL_{m,n}$
  the {\it Berezinian} of $A$ 
is given by the formula
$$\Ber A=\det (A_{11}-A_{12}A_{22}^{-1}A_{21})
\det A_{22}^{-1}=
$$
$$
=\det A_{11}
\det (A_{22}-A_{21}A_{11}^{-1}A_{12})^{-1}
$$
(see \cite {Ber}).
The main feature of the Berezinian 
is its mutiplicativity.

Denote
$$
T_l(V)=\underbrace {V\otimes \cdots\otimes V}_l,
$$
where $l=1,2,\ldots$.
Clearly, $T_l(V)$ is a $G$-module. 
By definition, for 
any ${\cal A}\in GL\, V$, $v_i\in V$,
 $l=1,2,\ldots$
we have 
\begin{equation}
{\cal A}(v_1 \cdots  v_l)=
{\cal A}(v_1) \cdots  {\cal A} (v_l).
\label{1.1}
\end{equation}
Action (\ref{1.1}) of $GL\, V$ on 
decomposable tensors
is extended to the whole space $T_l(V)$ by
 linearity and we see that
 $T_l(V)$ is a  $GL\,V$-module.

In what follows we consider the action 
 induced by $GL\;V$ on some submodules of  
$T_l(V)$. 
Nevertheless again the corresponding operator 
  is denoted by
$\cal A$. This cannot confuse us forasmuch as
 from
 context is clear, what of submodules 
is  considered.

Let $S_l$ be the symmetric group on the elements 
$ 1,2,\ldots ,l$.

By definition, for an arbitrary $\sigma\in S_l$,
$v_i\in V$, $i=1,2,\ldots,l$ put

\begin{equation}
\sigma (v_1\cdots v_l)=v_{\sigma^{-1}(1)}
\cdots v_{\sigma^{-1}(l)}.
\label{1.1a}    
\end{equation}
Then $T_l(V)$ is a left $S_l$-module.

To determine  irreducible $GL\,V$-submodules of 
$T_l(V)$  the special 
element of the group algebra $K[S_l]$,
the so-called  {\it  Young symmetrizer},
can be used.
Recall this notion.

Let $\lambda_1,\lambda_2,\ldots ,\lambda_k$ be
 positive integers
 such that $\lambda_1\geq \ldots \geq \lambda_k$
and $l=\lambda_1+\cdots +\lambda_k$. 
Then $\lambda=(\lambda_1,\ldots ,\lambda_k)$ is
{\it the partition} of $l$. The corresponding
 {\it Young diagram} $\lambda$ has $\lambda_i$
boxes in its $i$-th row. The {\it Young tableau} 
$T_\lambda$ is obtained from $\lambda$ by
filling its boxes with the numbers
 $1,2,\ldots ,l$ moving by the rows from left
to right and from top  to bottom.
For example, if $\lambda=(4,2,1)$, then $T_{\lambda}$
has the form
\begin{center}
\begin{picture}(65,50)
\put(0,48){\line(1,0){64}}
\put(0,32){\line(1,0){64}}
\put(0,16){\line(1,0){32}}
\put(0,0){\line(1,0){16}}
\put(0,0){\line(0,1){48}}
\put(16,0){\line(0,1){48}}
\put(32,16){\line(0,1){32}}
\put(48,32){\line(0,1){16}}
\put(64,32){\line(0,1){16}}
\put(5,36){$1$}
\put(21,36){$2$}
\put(37,36){$3$}
\put(52,36){$4$}
\put(5,20){$5$}
\put(21,20){$6$}
\put(5,4){$7$}
\end{picture}
\end{center}

Let $C(T_\lambda)$ ($R(T_\lambda$)) be the
 subgroup
 of $S_l$ that preserves the numbers of
 $T_\lambda$ within their columns
 (rows, respectively). Then 
$$
e_T=\sum_{
\begin{array}{c}
\sigma\in C(T),\\
\tau \in R(T)
\end{array}}
\sgn (\sigma)\tau \sigma
$$
is the {\it  Young  symmetrizer}
and  $W_\lambda=e_{T_\lambda}T_l(V)$
is the irreducible $GL\,V$-module.

Denote
$$
\lambda_h=(\underbrace {n+1,\ldots ,n+1}_m), 
$$
$$
\lambda_g=(\underbrace {n,\ldots ,n}_{m+1}).
$$
In other words, ${\lambda_h}$  is the
 $m\times (n+1)$ rectangle and  ${\lambda_g}$
is the $(m+1)\times n$ rectangle. In the sequel,
the $GL\,V$-modules $W_{\lambda_h}$,
 $W_{\lambda_g}$
play an important role.

{\bf Remark.}  One can say that 
 $W_{\lambda_g}$, $W_{\lambda_h}$ 
are particular cases of the so-called 
{\it external forms }
 (see \cite {Det}).

\vspace{0.5cm}

Denote
$$
\lambda_{SR}=(\underbrace {n,\ldots ,n}_m),
$$
$$
\lambda_{LR}=
(\underbrace {n+1,\ldots ,n+1}_{m+1}),
$$

$$
\lambda_{AR}=
(\underbrace {n+1,\ldots ,n+1}_m,n).
$$
According to terminology of Issaia L. Kantor
the diagrams ${\lambda_{SR}}$, 
${\lambda_{LR}}$, ${\lambda_{AR}}$
are called the {\it  small rectangle},
the {\it  large rectangle},
the {\it  almost rectangle} (respectively).

Evidently, to obtain ${\lambda_h}$
 (${\lambda_g}$) one must add the column of
  height  $m$ (the row of  length $n$)
to the small rectangle.
As mentioned above, in either case
 we obtain
a kind of rectangle.

An arbitrary tensor of the form
\begin{equation}
e_T v_1\cdots v_l, 
\label{1.2}    
\end{equation}
where $v_1,v_2,\ldots ,v_l\in V$, 
has the following properties:

i) if $i_1,i_2\in \{ 1,2,\ldots ,l\}$ place 
in the same column of $T$ 
then  tensor (\ref {1.2}) is skew-symmetric by
the elements $v_{i_1}$, $v_{i_2}$, that is, 
\begin{equation}
e_T v_1\cdots v_l=-e_T(i_1,i_2)v_1\cdots v_l; 
\label{1.3}    
\end{equation}

ii) let $i_1,\ldots ,i_{k-1}$ be  numbers 
that fill  a column of  length
 $(k-1)$ of $T$,
$i_k$ a number that belongs to a column of 
 length not more than $(k-1)$; then the {\it  
Jacobi
 identity} holds, that is, 
\begin{equation}
e_T v_1\cdots v_l=
\sum_{j=1}^{k-1}e_T(i_j,i_k)v_1\cdots v_l 
\label{1.4}    
\end{equation}
(see \cite {Form}).
Clearly, (\ref{1.3}) and (\ref{1.4}) 
are the applications of the corresponding properties
of the symmetrizer $e_T$.

Suppose $f$ is a tensor of the form  (\ref{1.2}).
We say that $(i,j)$ is an $f$-{\it box} if this box
 belongs to the tableau $T$ defining the tensor $f$.
By analogy the concepts of an $f$-{\it row} and an
$f$-{\it column} are introduced.

Let
 $\{e_1,\ldots,e_m\}$, $\{\varepsilon_1,\ldots,
\varepsilon_n\}$ be some bases of the $G$-modules
 $V_{0}$ 
and $V_{1}$ respectively.
Denote
$$
\Omega_0=\{e_1,\ldots,e_m\},\;\;
\Omega_1=\{\varepsilon_1,\ldots,
\varepsilon_n\},\;\;
\Omega=\Omega_0\cup \Omega_1.
$$

We say that an element $f\in T_l(V)$ is an 
{\it $\Omega$-tensor} if $f$ has the form (\ref{1.2}),
where $v_i\in \Omega$ for $i=1,2,\ldots,l$.

Recall that we use the unique way of  filling
 of Young diagrams (see above).
Thus there is a one-to-one correspondence
between   $\Omega$-tensors  and 
fillings of the 
corresponding
Young diagrams with the vectors $v_i\in \Omega$.

A filling of a rectangular Young diagram by vectors
 $v_i\in \Omega$ is
 {\it canonical} if every even basis vector 
$e_i\in \Omega_0$ cannot be placed
anywhere except   boxes of the 
$i$-th row and every odd basis vector
$\varepsilon_j\in \Omega_1$ cannot 
be placed anywhere except    boxes
of the $j$-th column. In other words, 
a filling of a rectangular diagram is canonical if a 
box $(i,j)$ that is within the small rectangle
 can be filled by $e_i$ or by $\varepsilon_j$
only. 

For example, if the dimension of $V$ is $2|1$,
$\Omega_0=\{e_1,e_2\}$, $\Omega_1=\{\varepsilon\}$,
then the canonical fillings of the small rectangle 
are the following:

\begin{center}
\begin{picture}(235,40)
\put(32,8){\line(1,0){16}}
\put(32,24){\line(1,0){16}}
\put(32,40){\line(1,0){16}}
\put(32,8){\line(0,1){32}}
\put(48,8){\line(0,1){32}}
\put(35,14){$e_2$}
\put(35,30){$e_1$}
\put(50,24){,}

\put(80,8){\line(1,0){16}}
\put(80,24){\line(1,0){16}}
\put(80,40){\line(1,0){16}}
\put(80,8){\line(0,1){32}}
\put(96,8){\line(0,1){32}}
\put(85,14){$\varepsilon$}
\put(83,30){$e_1$}
\put(98,24){,}

\put(128,8){\line(1,0){16}}
\put(128,24){\line(1,0){16}}
\put(128,40){\line(1,0){16}}
\put(128,8){\line(0,1){32}}
\put(144,8){\line(0,1){32}}
\put(133,30){$\varepsilon$}
\put(131,14){$e_2$}
\put(146,24){,}

\put(176,8){\line(1,0){16}}
\put(176,24){\line(1,0){16}}
\put(176,40){\line(1,0){16}}
\put(176,8){\line(0,1){32}}
\put(192,8){\line(0,1){32}}
\put(181,14){$\varepsilon$}
\put(181,30){$\varepsilon$}
\put(194,24){.}

\end{picture}
\end{center}

Suppose $\lambda$ is a rectangular Young diagram
filled with elements of $\Omega$ in a canonical way.
Then
the contents of boxes that do
 not belong to the small rectangle is defined
 uniquely. More precisely, if these boxes belong
to the $i$-th
 row, where $i\in \{1,2,\ldots,m\}$ 
($j$-th column, where $j\in \{1,2,\ldots,n\}$),
 then they are filled
 with $e_i$ (with $\varepsilon_j$ respectively).

In particular, the canonical fillings of the
 diagrams $\lambda_g$, $\lambda_h$ have the
 following properties:

i) the $(n+1)$-th column of the diagram
 $\lambda_h$ is filled with the vectors
$e_1,\ldots,e_m$ from the top down;

ii) the $(m+1)$-th row of the diagram 
$\lambda_g$ is filled with the vectors 
$\varepsilon_1,\ldots,
\varepsilon_n$ from left to right.

For example, if the dimension of $V$ is $1|2$,
$\Omega_0=\{e\}$, 
$\Omega_1=\{\varepsilon_1,\varepsilon_2\}$,
then the canonical fillings of the diagram 
$\lambda_g$ are the following:

\begin{center}
\begin{picture}(230,32)
\put(0,0){\line(1,0){32}}
\put(0,16){\line(1,0){32}}
\put(0,32){\line(1,0){32}}
\put(0,0){\line(0,1){32}}
\put(16,0){\line(0,1){32}}
\put(32,0){\line(0,1){32}}
\put(6,22){$e$}
\put(22,22){$e$}
\put(4,6){$\varepsilon_1$}
\put(20,6){$\varepsilon_2$}
\put(34,16){,}

\put(64,0){\line(1,0){32}}
\put(64,16){\line(1,0){32}}
\put(64,32){\line(1,0){32}}
\put(64,0){\line(0,1){32}}
\put(80,0){\line(0,1){32}}
\put(96,0){\line(0,1){32}}
\put(68,22){$\varepsilon_1$}
\put(86,22){$e$}
\put(68,6){$\varepsilon_1$}
\put(84,6){$\varepsilon_2$}
\put(98,16){,}

\put(128,0){\line(1,0){32}}
\put(128,16){\line(1,0){32}}
\put(128,32){\line(1,0){32}}
\put(128,0){\line(0,1){32}}
\put(144,0){\line(0,1){32}}
\put(160,0){\line(0,1){32}}
\put(134,22){$e$}
\put(148,22){$\varepsilon_2$}
\put(132,6){$\varepsilon_1$}
\put(148,6){$\varepsilon_2$}
\put(162,16){,}

\put(192,0){\line(1,0){32}}
\put(192,16){\line(1,0){32}}
\put(192,32){\line(1,0){32}}
\put(192,0){\line(0,1){32}}
\put(208,0){\line(0,1){32}}
\put(224,0){\line(0,1){32}}
\put(196,22){$\varepsilon_1$}
\put(212,22){$\varepsilon_2$}
\put(196,6){$\varepsilon_1$}
\put(212,6){$\varepsilon_2$}
\put(226,16){.}

\end{picture}
\end{center}

 Evidently, the concept
of  canonical filling is defined for 
 diagrams
 that do not contain the large rectangle.
From the other hand, if a diagram  
 contains the large rectangle,
then tensor (\ref{1.2}) is equal to zero 
(see, for example, \cite{Det}).

{\bf Remark.}  The concept of canonical filling
is introduced in  \cite{Det} and can be applied
not only to rectangular diagrams.

\vspace{0.5cm}

We say that an $\Omega$-tensor $f$ is
{\it canonical} if the filling of the
 corresponding diagram is canonical.

By $\Lambda_g$, $\Lambda_h$ denote the sets 
of all canonical tensors for the
 diagrams $\lambda_g$,
$\lambda_h$. Note that there is one-to-one
 correspondence
between the sets $\Lambda_g$, $\Lambda_h$: 
the corresponding elements have the same
 filling  within the small rectangle.
Elements of the sets $\Lambda_g$, $\Lambda_h$
are denoted by $g_i$ and $h_j$ respectively,
where $i,j=1,2,\ldots ,k_\Lambda$, 
$k_\Lambda=2^{mn}$.
Moreover we suppose that elements of 
$\Lambda_g$, $\Lambda_h$ with the same numbers
are corresponding, that is, for every  
$i\in \{ 1,2,\ldots,k_\Lambda\}$
 tensors $g_i$, $h_i$ have the same  
 filling of the small rectangle.

{\bf Theorem 1.1.} 
 {\em The sets $\Lambda_g$,
$\Lambda_h$ are the bases of the $G$-modules 
$W_{\lambda_g}$, $W_{\lambda_h}$. 
}

The proof of this theorem is found in  \cite {Det}.

\section{{\rm The isomorphism of  $W_{\lambda_h}$
and $ b \cdot W_{\lambda_g}$ as the 
$GL\, V$-modules }}
\setcounter{equation}{0}
\indent

 Let $b$ be a formal element such that
 for any
$
{\cal A} \in GL\, V
$
we have
\begin{equation}
 {\cal A}(b)=\Ber A\cdot b,
\label{2.1}    
\end{equation} 
where $A$ is the matrix of $\cal A$
in some base of $V$.
Recall that to change a base of $V$ 
we use even operators only (see \S 1). 
Then taking into account the multiplicativity 
of $\Ber$, one can see that 
the concept of the element
$b$ is well defined.

Clearly, we assume the element $b$ is
homogeneous. At the same time now we do not
 define if this element even or odd.

Because of the multiplicativity of 
$\Ber$ formula (\ref {2.1}) defines the 
one-dimensional
 representation of the group $GL\, V$ 
on the  $G$-module generated by
 the element $b$.

Let $f$ be an $\Omega$-tensor.
By $\kappa(f)$ denote the cardinality of the 
automorphism group of $f$-tableau by odd elements. 
Clearly, if an $f$-column does not contain
any odd elements, then the corresponding factor
in $\kappa(f)$ is equal to $1$. 

By $\rho (f)$ denote the
 number of $\Omega_1$-elements belonging to the
 small rectangle of $f$. Let $\alpha$ be the map
from $\Lambda_h$ to the rational numbers such that
\begin{equation}
 \alpha(h_i)=\frac{(-1)^{n\rho(h_i)}\kappa(h_i)}
{\kappa(g_i)}.
\label{2.2}    
\end{equation}

By definition, put

\begin{equation}
h_i^{\prime}=\frac{h_i}{\alpha(h_i)},
\label{2.3}    
\end{equation}
where $i=1,2,\ldots,k_{\Lambda}$.

The following result is the heart of the 
present paper.

{\bf Theorem 2.1.} 
 {\em The mapping
\begin{equation}
\varphi : h_i^{\prime}\mapsto b\cdot g_i
\label{2.4}    
\end{equation}
determines the  isomorphism of the $GL\; V$-modules
$W_{\lambda_h}$ and $b\cdot W_{\lambda_g}$, that is,
for any ${\cal A}\in {GL\; V}$ the following 
equality holds
\begin{equation}
{\cal A}\varphi=\varphi{\cal A}.
\label{2.5}    
\end{equation}
 }

In the sequel we assume that  
$A\in GL_{m,n}$ is the matrix
of ${\cal A}\in {GL\; V}$  in the base 
$ e_1,\ldots,e_m,\varepsilon_1,\ldots,
\varepsilon_n$.

The following proposition    reduces the proof
of Theorem 2.1 to some easier particular cases of
the general situation.

{\bf Proposition 2.1.}  {\em  To prove Theorem 2.1
 it  suffices to check (\ref{2.5}) for the cases
 when the matrix $A$ of ${\cal A}\in {GL\; V}$ 
belongs to one of the following classes of matrices:

\begin{equation}
 E_{m+n}+e_{ij}\eta_{ij},
\label{2.6}    
\end{equation}
where $E_{m+n}$ is the unit matrix of order $m+n$,
$e_{ij}$ is the matrix unit, $\eta_{ij}\in G_1$,
$m+1\leq i\leq m+n$, $1\leq j\leq m$;

\begin{equation}
 E_{m+n}+e_{ij}\xi_{ij},
\label{2.7}    
\end{equation}

where $\xi_{ij}\in G_1$, $1\leq i\leq m$,
$m+1\leq j\leq m+n$;

\begin{equation}
\diag (x_1,\ldots,x_m,y_1,\ldots,y_n),
\label{2.8}    
\end{equation}
where $x_i,y_j\in G_0$;

\begin{equation}
 E_{m+n}+e_{ij}x_{ij},
\label{2.9}    
\end{equation}
where $x_{ij}\in G_0$, $1\leq i,j\leq m$, $i\ne j$;

\begin{equation}
 E_{m+n}+e_{ij}y_{ij},
\label{2.10}    
\end{equation}
where $y_{ij}\in G_0$,
 $m+1\leq i,j\leq m+n$, $i\ne j$.}

The proof of this proposition is based on the 
following auxiliary result.

{\bf Lemma 2.1.} {\em Any matrix 
$A=\left( \begin {array}{l} 
A_{11}\;A_{12}\\A_{21}\;A_{22}
\end{array} \right)\in GL_{m,n}$ can be represented
in the form

\begin{equation}
A=\left( \begin {array}{cc} 
E_m &0\\B_{21}&E_n
\end{array} \right)
\left( \begin {array}{cc} 
B_{11}&0\\0&B_{22}
\end{array} \right)
\left( \begin {array}{cc} 
E_m&B_{12}\\0&E_n
\end{array} \right),
\label{2.11}    
\end{equation}

where $B_{11}=A_{11}$, $B_{12}=A_{11}^{-1}A_{12}$,
$B_{21}=A_{21}A_{11}^{-1}$,
$B_{22}=A_{22}-A_{21}A_{11}^{-1}A_{12}$,

or in the form

\begin{equation}
A=\left( \begin {array}{cc} 
E_m&C_{12}\\0&E_n
\end{array} \right)
\left( \begin {array}{cc} 
C_{11}&0\\0&C_{22}
\end{array} \right)
\left( \begin {array}{cc} 
E_m&0\\C_{21}&E_n
\end{array} \right),
\label{2.12}    
\end{equation}
where $C_{22}=A_{22}$, $C_{12}=A_{12}A_{22}^{-1}$,
$C_{21}=A_{22}^{-1}A_{21}$,
$C_{11}=A_{11}-A_{12}A_{22}^{-1}A_{21}$.} 

{\sc Proof.} Since $A\in GL_{m,n}$, we see that 
$A_{11}^{-1}$, $A_{22}^{-1}$ are exist
 (see \cite {Ber}). Equalities (\ref{2.11}),
(\ref{2.12}) can be checked directly.

This completes the proof.

{\sc Proof of Proposition 2.1.} With the isomorphism 
$GL\; V\simeq GL_{m,n}$ a composition of linear operators
corresponds to the composition of matrices 
(clearly, the inverse statement is correct).
Hence it follows from 
Lemma 2.1 that it  suffices to check (\ref {2.5})
for  ${\cal A}\in GL\; V$ such that the matrix $A$
of $\cal A$ has either the form

\begin{equation}
\left( 
\begin {array}{cc}
 E_m&A_{12}
\\
0&E_n
\end{array} 
\right),
\label{2.13}    
\end{equation}
or
\begin{equation}
\left( 
\begin {array}{cc}
E_m&0
\\
A_{21}&E_n
\end{array} 
\right),
\label{2.14}    
\end{equation}
or
\begin{equation}
\left( 
\begin {array}{cc}
A_{11}&0
\\
0&A_{22}
\end{array} 
\right),
\label{2.15}    
\end{equation}
where (as above) $A_{11}$, $A_{22}$ are square 
$G_0$-matrices of orders $m$ and $n$ respectively,
$A_{12}$, $A_{21}$ are $G_1$-matrices.

If $i$, $j$  are  integers such that $1\leq i\leq m$,
$m+1\leq j\leq m+n$,  then $e_{ij}^2=0$.
Consequently an arbitrary matrix 
(\ref{2.13}) can be represented as a product 
of matrices  (\ref {2.7}).
Similarly any matrix  (\ref{2.14})
can be represented as a product of matrices 
 (\ref{2.6}).

By $\cal D$ denote the group generated by all matrices
of the  form (\ref{2.8}) - (\ref{2.10}).
Also by $\bar {\cal D}$ denote the group of all matrices 
of the  form (\ref {2.15}).
Just as for  the regular subgroup in the classic case
 the set  $\cal D$
is open and dense in $\bar {\cal D}$.

Thus we see that Proposition 2.1 is proved.

Applying Proposition 2.1 we use the coordinate
 form of equality (\ref{2.5}).
Let us obtain this form. 
With (\ref{2.1}) and (\ref {2.4}) we have 
$$
\Ber A\cdot b\cdot {\cal A}(g_i)=
{\cal A}(b){\cal A}(g_i)=
$$
$$
={\cal A}(bg_i)={\cal A}\varphi(h_i^{\prime}).
$$
Hence (\ref{2.5}) is equivalent to
\begin{equation}
\varphi({\cal A}(h_i^{\prime}))=\Ber A \cdot
b\cdot{\cal A}(g_i).
\label{2.16}    
\end{equation}

The matrix elements
of the action induced by $\cal A$ on the spaces
$W_{\lambda_h}$, $W_{\lambda_g}$ 
and written in the bases $\{h_i^{\prime}\}$ and
$\{g_j\}$
we denote by
 $a_{ij}^{\prime}$ and $a_{ij}^{\prime\prime}$ 
(respectively), that is,

\begin{equation}
{\cal A}(h_j^{\prime})=\sum_{i=1}^{k_\Lambda}
h_i^{\prime}a_{ij}^{\prime},
\label{2.17}    
\end{equation}

\begin{equation}
{\cal A}(g_j)=\sum_{i=1}^{k_\Lambda}
g_i a_{i j}^{\prime \prime},
\label{2.18}    
\end{equation}
where $1\leq j\leq k_\Lambda$.

Using (\ref{2.17}), (\ref{2.18}) rewrite 
condition (\ref{2.16}) in the form
$$
\varphi
\left(
\sum_{i=1}^{k_\Lambda}h_i^{\prime}
a_{i j}^{\prime}
\right)
=
\Ber A\cdot b\cdot 
\sum_{i=1}^{k_\Lambda}g_i a_{ij}^{\prime\prime}.
$$

Hence by definition of $\varphi$ we have

\begin{equation}
b\cdot\sum_{i=1}^{k_\Lambda}
g_ia_{ij}^{\prime}=
\Ber A\cdot b\cdot
\sum_{i=1}^{k_\Lambda}
g_ia_{ij}^{\prime\prime}.
\label{2.20}    
\end{equation}

Since the elements $g_i$ are $G$-linearly
 independent,
we see that (\ref{2.20}) is equivalent to
\begin{equation}
a_{ij}^{\prime}=\Ber A\cdot a_{ij}^{\prime\prime},
\label{2.21}    
\end{equation}
where $1\leq i,j\leq k_\Lambda$. Equality
 (\ref{2.21})
is the  coordinate form  of (\ref{2.5}).

For an arbitrary ${\cal A}\in GL\; V$ by 
$A_{h^{\prime}}$, $A_g$ denote the matrices 
of the action induced by ${\cal A}$  on the spaces
$W_{\lambda_h}$, $W_{\lambda_g}$ and 
written in the bases $\{ h_i^{\prime}\}$, 
$\{ g_i\}$, that is, 
$A_{h^{\prime}}=(a_{ij}^{\prime})$, 
$A_g=(a_{ij}^{\prime\prime})$.

{\bf Lemma 2.2.} {\em  If the matrix $A$ of
${\cal A}\in GL\; V$ has the form  (\ref{2.6}),
then 
$$
A_{h^{\prime}}=A_g.
$$}

For any matrix $A$ of the form (\ref{2.6})
we have $\Ber A=1$. Hence, according to
 Proposition 2.1, Lemma 2.2 is a part of
 the proof of Theorem 2.1 (see (\ref{2.21})).

{\sc Proof of Lemma 2.2.}  By assumption,
 there exist integers $i\in \{ 1,\ldots,n\}$,
$j\in \{1,\ldots,m\}$ such that 

$$
{\cal A}(e_j)=e_j+\varepsilon_i\eta_{ij},
$$
where $\eta_{ij}\in G_1$, and ${\cal A}(w)=w$
for $w\in \Omega\backslash \{e_j\}$.

Let $h_t$, $g_t$ be corresponding elements
of the bases $\{ h_i\}$, $\{ g_i \}$, where
$1\leq t\leq k_\Lambda$, 
$Q_{g_t}=\{ j_1,\ldots ,j_l\}$
$(Q_{h_t})$ the set of the numbers of $g_t$-columns
($h_t$-columns, respectively) that contain $e_j$.
In other words, $r\in Q_{g_t}$ $(Q_{h_t})$
iff the $g_t$-box ($h_t$-box, resp.) $(j,r)$
is filled with $e_j$.
Evidently,  $Q_{h_t}=Q_{g_t}\cup \{ n+1\}$.
Then we have
\begin{equation}
{\cal A}(h_t)=h_t+
\sum_{r\in Q_{h_t}}{\hat h}_r\gamma_{j,r}\eta_{ij},
\label{2.23}    
\end{equation}
where  ${\hat h}_r$ is obtained from $h_t$
by the substitution of $\varepsilon_i$ for
 $e_j$ placed  in the $h_t$-box
$(j,r)$;

\begin{equation}
\gamma_{j,r}=(-1)^{n_{j,r}},
\label{2.24}    
\end{equation}
where $n_{j,r}$ is the number of odd base vectors in
 the $h_t$-tableau 
within the interval from the box
$(j,r+1)$ to the box $(m,n+1)$
(recall that we move by the tableau from left
to right and from top to bottom by  rows).
Clearly, the tensors ${\hat h}_r$ depend on 
$t$, $i$, $j$. Nonetheless, we omit these indices in
 the notation  to simplify the last one.
 In general, the tensors ${\hat h}_r$
are not canonical.

Suppose $i\neq Q_{h_t}$, that is, the 
$h_t$-box $(j,i)$  is filled with $\varepsilon_i$.

Then any one of $\{{\hat h}_r|r\in Q_{h_t}\}$
is not canonical. To express ${\hat h}_r$ in terms of
 canonical tensors we use the Jacobi identity 
for the $i$-th column and for the box $(j,r)$ of 
the ${\hat h}_r$-tableau. Also we use (\ref{1.3}) 
to transpose $e_p$ and $\varepsilon_r$ that are 
placed in the boxes $(j,r)$ and $(p,r)$ 
after applying (\ref{1.4}).

Thus we obtain

\begin{equation}
{\hat h}_r=-q{\hat h}_r+\sum_{p\in P_r}
h_{r,p}\beta_{r,p},
\label{2.25}    
\end{equation}
where  $q$ is the number of  $\varepsilon_i$
in the $i$-th $h_t$-column  
(or, what is the same, in the
$i$-th ${\hat h}_r$-column);
$P_r\subseteq \{1,\ldots,m\}$ is the set of integers
  such that $z\in P_r$ iff $e_z$ belongs to the 
$i$-th ${\hat h}_r$-column and $e_z$ does not belong
to the $r$-th ${\hat h}_r$-column.
The tensor $h_{r,p}$ is obtained from ${\hat h}_r$
in the following two steps:

i)  the elements $\varepsilon_i$ and
$e_p$ disposed in the ${\hat h}_r$-boxes $(j,r)$ and
$(p,i)$ (respectively) change places;

ii)  the elements $e_p$ and $\varepsilon_r$
disposed in the boxes $(j,r)$ and $(p,r)$ (resp.)
change places (see Fig. 1).

$$
\beta_{r,p}=(-1)^{{\breve m}_{r,p}+1},
$$
where ${\breve m}_{r,p}$ is the number of pairs
of odd elements that change their order under 
 transpositions i) - ii); the additional
term $(+1)$ in the exponent appears in the 
application of (\ref{1.3}). 
 
It is easy to see that the 
transformation of ${\hat h}_r$-tableau to  
 $h_{r,p}$-tableau also can be realised in
 the following two steps:
first the elements $\varepsilon_i$ and 
$\varepsilon_r$ disposed in the boxes $(j,r)$ and
$(p,r)$ (respectively) change places;
then the elements $\varepsilon_i$ and
$e_p$ disposed in the boxes $(p,r)$ and $(p,i)$
(resp.) change places (see Fig. 2).

\begin{center}
\begin{picture}(320,140)

\put(20,10){\line(1,0){120}}
\put(20,30){\line(1,0){120}}
\put(20,50){\line(1,0){120}}
\put(20,80){\line(1,0){120}}
\put(20,100){\line(1,0){120}}
\put(20,120){\line(1,0){120}}
\put(20,10){\line(0,1){110}}
\put(50,10){\line(0,1){110}}
\put(70,10){\line(0,1){110}}
\put(100,10){\line(0,1){110}}
\put(120,10){\line(0,1){110}}
\put(140,10){\line(0,1){110}}
\put(10,37){$p$}
\put(10,87){$j$}
\put(56,37){$\varepsilon_r$}
\put(56,87){$\varepsilon_i$}
\put(57,125){$r$}
\put(106,37){$e_p$}
\put(106,87){$\varepsilon_i$}
\put(107,125){$i$}
\put(53,45){\vector(0,1){42}}
\put(67,90){\vector(0,-1){46}}
\put(105,40){\vector(-3,4){38}}
\put(67,95){\vector(3,-4){39}}
\put(49,79){\line(0,1){22}}
\put(49,79){\line(1,0){22}}
\put(71,79){\line(0,1){22}}
\put(49,101){\line(1,0){22}}
\put(99,9){\line(1,0){22}}
\put(99,9){\line(0,1){112}}
\put(121,9){\line(0,1){112}}
\put(99,121){\line(1,0){22}}
\put(66,0){Fig. 1}

\put(180,10){\line(1,0){120}}
\put(180,30){\line(1,0){120}}
\put(180,50){\line(1,0){120}}
\put(180,80){\line(1,0){120}}
\put(180,100){\line(1,0){120}}
\put(180,120){\line(1,0){120}}
\put(180,10){\line(0,1){110}}
\put(210,10){\line(0,1){110}}
\put(230,10){\line(0,1){110}}
\put(280,10){\line(0,1){110}}
\put(280,10){\line(0,1){110}}
\put(300,10){\line(0,1){110}}
\put(170,37){$p$}
\put(170,87){$j$}
\put(216,37){$\varepsilon_r$}
\put(216,87){$\varepsilon_i$}
\put(217,125){$r$}
\put(266,37){$e_p$}
\put(266,87){$\varepsilon_i$}
\put(267,125){$i$}
\put(213,45){\vector(0,1){42}}
\put(227,87){\vector(0,-1){43}}
\put(227,44){\vector(1,0){39}}
\put(265,33){\vector(-1,0){40}}
\put(209,79){\line(0,1){22}}
\put(209,79){\line(1,0){22}}
\put(231,79){\line(0,1){22}}
\put(209,101){\line(1,0){22}}
\put(259,9){\line(1,0){22}}
\put(259,9){\line(0,1){112}}
\put(281,9){\line(0,1){112}}
\put(259,121){\line(1,0){22}}
\put(226,0){Fig. 2}
\end{picture}
\end{center}

Then we see that 
$$
\beta_{r,p}=(-1)^{m_{r,p}},
$$
where $m_{r,p}$ is the number of  odd
base vectors of the $p$-th row  that are placed
strictly between the $r$-th and the $i$-th columns.

Now it follows from (\ref{2.25}) that

\begin{equation}
{\hat h}_r=\frac{1}{q+1}\sum_{p\in P_r}
h_{r,p}\beta_{r,p}.
\label{2.27}    
\end{equation}
Notice that the tensors $h_{r,p}$ in (\ref{2.27}) 
are canonical.
Note also that  $P_{n+1}=\emptyset$.
Hence using the Jacobi identity we get
${\hat h}_{n+1}=0$.

 Combining (\ref{2.27}) with (\ref{2.23}) we get
\begin{equation}
{\cal A}(h_t)=
h_t+\frac{1}{q+1} 
\sum_{r\in Q}\gamma_{j,r}
\sum_{p\in P_r}h_{r,p}\beta_{r,p}
\eta_{ij},
\label{2.28}    
\end{equation}
where $Q=Q_{h_t}\backslash \{n+1\}$.

By assumption, $i\notin Q$. We claim that if  
$(r_1,p_1)\neq(r_2,p_2)$, 
then $h_{r_1,p_1}\neq h_{r_2,p_2}$.
Indeed, if $r_1\neq r_2$ then the sets of columns 
containing $e_j$ are distinct for $h_{r_1,p_1}$
and $h_{r_2,p_2}$. Further, an element $e_p$, 
where $p\in P_r$, belongs to the $r$-th 
$h_{r,z}$-column iff $z=p$. Hence, 
$h_{r,p_1}\neq h_{r,p_2}$ for $p_1\neq p_2$.

In terms of the ``primed'' tensors (see (\ref{2.3}))
 equality (\ref{2.28}) takes the form 
\begin{equation}
{\cal A}(h_t^{\prime})=
h_t^{\prime}+\frac{1}{(q+1)\alpha_t}
\sum_{r\in Q}\gamma_{j,r}
\sum_{p\in P_r}
h_{r,p}^{\prime}\beta_{r,p}\alpha_{r,p}
\eta_{ij},
\label{2.29}    
\end{equation}
where $\alpha_t=\alpha(h_t)$,
$\alpha_{r,p}=\alpha(h_{r,p})$ 
(see (\ref{2.2})).

Similarly, still assuming that $i\notin Q$
we get
\begin{equation}
{\cal A}(g_t)=
g_t+\sum_{r\in Q}
{\hat g}_r{\bar \gamma_{j,r}}\eta_{ij},
\label{2.30}    
\end{equation}
where ${\hat g}_r$ is obtained from
 $g_t$ by replacing the vector $e_j$ 
in the box $(j,r)$
to $\varepsilon_i$;

\begin{equation}
{\bar \gamma}_{j,r}=(-1)^{{\bar n}_{j,r}},
\label{2.30.1}    
\end{equation}
where ${\bar n}_{j,r}$ is the number of
 odd base vectors that are placed after
the $g_t$-box $(j,r)$ when we move by the tableau
from left to right and from the top down.
Since the tableaux of $h_t$ and $g_t$
have the same filling within the small rectangle
and the $(m+1)$-th $g_t$-row of length $n$
is filled with the odd base vectors 
$\varepsilon_1,\ldots,\varepsilon_n$,
we have 

\begin{equation}
{\bar n}_{j,r}=n_{j,r}+n.
\label{2.30.2}    
\end{equation}
With (\ref{2.30.1}), (\ref{2.30.2}), (\ref{2.24})
equality (\ref{2.30}) takes the form

$$
{\cal A}(g_t)=g_t+
\sum_{r\in Q}
{\hat g}_r\gamma_{j,r}(-1)^n\eta_{ij},
$$
where $Q=Q_{g_t}=Q_{h_t}\backslash \{n+1\}$.

The tensors ${\hat g}_r$ are not canonical
because the element $\varepsilon_i$ is placed not in
 its ``native'' $i$-th column but in the $r$-th column
($r\neq i$).

To express ${\hat g}_r$ in terms of canonical tensors
we use again the Jacobi identity.
By analogy with (\ref{2.25}), writing this identity
for the $i$-th column and the box $(j,r)$ we get

\begin{equation}
{\hat g}_r=-(q+1){\hat g}_r
+\sum_{p\in P_r}
g_{r,p}{\bar\beta}_{r,p},
\label{2.31}    
\end{equation}
where $g_{r,p}$ is obtained from ${\hat g}_r$
in  the same way as  $h_{r,p}$
is obtained from ${\hat h}_r$ (see above);

$$
{\bar\beta}_{r,p}=(-1)^{{\bar m}_{r,p}},
$$
where ${\bar m}_{r,p}$ is defined by analogy with
$m_{r,p}$ (see above) but now for ${\hat g}_r$.
Since the tableaux of ${\hat h}_r$ and
${\hat g}_r$ have the same filling within the small
 rectangle, we have ${\bar m}_{r,p}=m_{r,p}$.
Hence,
$$
{\bar\beta}_{r,p}=\beta_{r,p}.
$$

From (\ref{2.31}) it follows that

\begin{equation}
{\hat g}_r=
\frac{1}{q+2}
\sum_{p\in P_r}
g_{r,p}\beta_{r,p}.
\label{2.32}    
\end{equation}

With (\ref{2.32}) equality (\ref{2.30}) 
takes the form  

\begin{equation}
{\cal A}(g_t)=g_t+
\frac{(-1)^n}{q+2}
\sum_{r\in Q}
\gamma_{j,r}
\sum_{p\in P_r}
g_{r,p}\beta_{r,p}\eta_{ij}.
\label{2.33}    
\end{equation}

We see  that the canonical tensors $h_{r,p}$ and
$g_{r,p}$ are corresponding in that these have
the same filling of the small rectangle.

When (\ref{2.29}) is compared with (\ref{2.33}),
it is apparent that
to conclude the proof it remains to check
the equality

\begin{equation}
\frac{\alpha_{r,p}}{\alpha_t}=
\frac{(-1)^n(q+1)}{q+2}.
\label{2.34}    
\end{equation}

With (\ref{2.3}) the last equality 
 is equivalent to
\begin{equation}
\frac{(-1)^{n\rho(h_{r,p})}
\kappa(h_{r,p})\kappa(g_t)}
{\kappa(g_{r,p})
(-1)^{n\rho(h_t)}
\kappa(h_t)}=
\frac{(-1)^n(q+1)}{q+2}.
\label{2.35}    
\end{equation}

First note that $\rho(h_{r,p})-\rho(h_t)=1$. 

Further, for an arbitrary $z\in \{1,\ldots,n\}$,
$z\neq i $ the number of odd base vectors placed
 in  the $z$-th  column is the same for $h_t$
and $h_{r,p}$ (also for $g_t$ and $g_{r,p}$).
The numbers of odd base vectors in the $i$-th column
 of $h_t$, $h_{r,p}$, $g_t$, $g_{r,p}$ are equal 
 to $q$, $q+1$, $q+1$, $q+2$ respectively.
Then we have
$$
\frac{\kappa(h_{r,p})\kappa(g_t)}
{\kappa(g_{r,p})\kappa(h_t)}=
\frac{(q+1)!(q+1)!}{(q+2)!q!}=
\frac{q+1}{q+2}.
$$

Hence (\ref{2.35}) is proved.

Thus Lemma 2.2 is proved for the case when
$i\notin Q$.

Now suppose $i\in Q$, that is, the $h_t$-box $(j,i)$
 is filled with $e_j$.
Still we have (\ref{2.23}). The tensors 
${\hat h}_r$ are not
canonical except when $r=i$.

The tensor ${\hat h}_{n+1}$ is not canonical 
because the element $\varepsilon_i$
is in the box $(j,n+1)$, but
the $(n+1)$-th column must be filled with even 
base vectors only.
To express  ${\hat h}_{n+1}$ in terms of canonical
 tensors we use the Jacobi identity for the 
$i$-th column and the box $(j,n+1)$.
 In such a manner we get

\begin{equation}
{\hat h}_{n+1}=-q{\hat h}_{n+1}
+{\bar h}\beta_{n+1,j},
\label{2.36}    
\end{equation}
where as above $q$ is the number of $\varepsilon_i$
placed in the $i$-th $h_t$-column;  
$\bar h$ is obtained from ${\hat h}_{n+1}$ by the
transposition of $e_j$ and $\varepsilon_i$
placed in the ${\hat h}_{n+1}$-boxes $(j,i)$ 
and $(j,n+1)$
(respectively).

Evidently, ${\bar h}={\hat h}_i$.
Also,
$$
\beta_{n+1,j}=(-1)^{m_{n+1,j}},
$$
where $m_{n+1,j}$ is the number of odd base vectors
placed in the $j$-th ${\hat h}_{n+1}$-row strictly 
between $i$-th 
and $(n+1)$-th columns.

Now it follows from (\ref{2.36}) that
\begin{equation}
{\hat h}_{n+1}=\frac{\beta_{n+1,j}}{q+1}
{\hat h}_i.
\label{2.38}    
\end{equation}

Consider an arbitrary tensor ${\hat h}_r$,
where $r\in Q\backslash \{ i \}$.
We claim that  $j\in P_r$, where $P_r$ is the set defined above.
Indeed, by assumption the ${\hat h}_r$-box $(j,i)$ 
  is filled with $e_j$,
also $\varepsilon_i$ is in the box $(j,r)$
and $e_j$ can be placed only in the $j$-th row
in the tableau of ${\hat h}_r$. 

By applying the Jacobi identity to the $i$-th column
and the box $(j,r)$ of the tableau of ${\hat h}_r$,
we get

\begin{equation}
{\hat h}_r=-q{\hat h}_r+
\sum_{p\in P_r\backslash\{ j \}}
h_{r,p}\beta_{r,p}+
{\hat h}_i\beta_{r,j},
\label{2.39}    
\end{equation} 
where $h_{r,p}$, $\beta_{r,p}$
are defined above.

Recall that ${\hat h}_r$ is canonical iff $r=i$.
By this reason (\ref{2.23}) is conveniently
rewritten in the form:

\begin{equation}
{\cal A}(h_t)=h_t+
\sum_{r\in Q\backslash \{i\}}
{\hat h}_r\gamma_{j,r}\eta_{ij}+
{\hat h}_i\gamma_{j,i}\eta_{ij}+
{\hat h}_{n+1}\gamma_{j,n+1}\eta_{ij}.
\label{2.40}    
\end{equation}

With (\ref{2.38}) and
(\ref{2.39}) equality (\ref{2.40}) takes the form:

$$
{\cal A}(h_t)=h_t+
\frac{1}{q+1}
\sum_{r\in Q\backslash \{ i \}}
\gamma_{j,r}
\sum_{p\in P_r\backslash\{ j \}}
h_{r,p}\beta_{r,p}
\eta_{ij}+
$$
\begin{equation}
\sum_{r\in Q\backslash \{ i \}}
\frac{\gamma_{j,r}}{q+1}{\hat h}_i
\beta_{r,j}\eta_{ij}+
{\hat h}_i\gamma_{j,i}\eta_{ij}+
\frac{\beta_{n+1,j}}{q+1}{\hat h}_i
\gamma_{j,n+1}\eta_{ij}.
\label{2.41}    
\end{equation}
Clearly (\ref{2.41}) is equivalent to
$$
{\cal A}(h_t)=h_t+
\frac{1}{q+1}
\sum_{r\in Q\backslash \{ i \}}
\gamma_{j,r}
\sum_{p\in P\backslash \{ j \}}
h_{r,p}\beta_{r,p}
\eta_{ij}+
$$

\begin{equation}
{\hat h}_i
\left(
\frac{1}{q+1}
\sum_{r\in Q\backslash \{ i \}}
\gamma_{j,r}\beta_{r,j}+
\gamma_{j,i}+
\frac{\beta_{n+1,j}}{q+1}
\gamma_{j,n+1}
\right)
\eta_{ij}.
\label{2.42}    
\end{equation}

Notice that $h_{r,p}\neq{\hat h}_i$,
where $r\in Q\backslash\{i\}$,
$p\in P_r\backslash\{j\}$.

We have

\begin{equation}
\gamma_{j,r}\beta_{r,j}=\gamma_{j,i},\;\;
\gamma_{j,n+1}\beta_{n+1,j}=\gamma_{j,i}.
\label{2.42-1}    
\end{equation}
Then the coefficient at ${\hat h}_i$ in
(\ref{2.42}) is equal to

$$
\frac{\gamma_{j,i}(l+q+2)}{q+1},
$$
where $l=|Q|$. Therefore equality
 (\ref{2.42}) takes the form

$$
{\cal A}(h_t)=h_t+
\frac{1}{q+1}
\sum_{r\in Q\backslash \{ i \}}
\gamma_{j,r}
\sum_{p\in P\backslash \{ j \}}
h_{r,p}\beta_{r,p}
\eta_{ij}
+\frac{l+q+2}{q+1}
{\hat h}_i\gamma_{j,i}\eta_{ij}.
$$

Going to the ``primed'' tensors (see (\ref{2.3}))
we obtain

$$
{\cal A}(h_t^{\prime})=h_t^{\prime}+
\frac{1}{(q+1)\alpha_t}
\sum_{r\in Q\backslash \{ i \}}
\gamma_{j,r}
\sum_{p\in P\backslash \{ i \}}
h_{r,p}^{\prime}\beta_{r,p}\alpha_{r,p}
\eta_{ij}+
$$

\begin{equation}
\frac{l+q+2}{(q+1)\alpha_t}
{\hat h}_i^{\prime}\gamma_{j,i}
{\hat \alpha}_i
\eta_{ij},
\label{2.45}    
\end{equation}
where ${\hat \alpha}_i=\alpha({\hat h}_i)$.

In its turn by analogy with (\ref{2.30}) we have

\begin{equation}
{\cal A}(g_t)=g_t+
\sum_{r\in Q\backslash \{ i \}}
{\hat g}_r
\gamma_{j,r}
\eta_{ij}(-1)^n+
{\hat g}_i\gamma_{j,i}\eta_{ij}(-1)^n,
\label{2.46}    
\end{equation}
where ${\hat g}_r$,  $\gamma_{j,r}$,
$\gamma_{j,i}$ are defined as above for the case when 
$i\notin Q$. An important point is that the tensors
${\hat g}_r$ are not canonical for 
$r\in Q\backslash \{ i \}$ and ${\hat g}_i$ is 
canonical.

To express  an arbitrary tensor ${\hat g}_r$, where 
$r\in Q\backslash \{ i \}$, in terms of canonical
 tensors we use the Jacobi identity for the $i$-th 
column and the box $(j,r)$. As it is mentioned above
the condition $i\in Q$ implies $j\in P_r$.
Then by analogy with (\ref{2.39}), (\ref{2.31})
we get

$$
{\hat g}_r=-(q+1){\hat g}_r+
\sum_{p\in P_r\backslash\{j\}}
g_{r,p} \beta_{r,p} +
{\hat g}_i\beta_{r,j},
$$

From the last equality it follows that

$$
{\hat g}_r=\frac{1}{q+2}
\sum_{p\in P_r\backslash\{j\}}
g_{r,p} \beta_{r,p} +
\frac{1}{q+2}
{\hat g}_i\beta_{r,j}.
$$

With the last expression for ${\hat g}_r$ equality
(\ref{2.46}) becomes:

$$
{\cal A}(g_t)=
g_t+
\frac{(-1)^n}{q+2}
\sum_{r\in Q\backslash\{i\}}
\gamma_{j,r}
\sum_{p\in P\backslash\{j\}}
g_{r,p}\beta_{r,p}\eta_{ij}+
$$

\begin{equation}
{\hat g}_i(-1)^n
\left(
\frac{1}{q+2}
\sum_{r\in Q\backslash\{i\}}
\gamma_{j,r}\beta_{r,j}
+\gamma_{j,i}
\right)
\eta_{ij}.
\label{2.49}    
\end{equation}

Using (\ref{2.42-1}) we obtain  

\begin{equation}
\frac{1}{q+2}
\sum_{r\in Q\backslash\{i\}}
\gamma_{j,r}\beta_{r,j}+
\gamma_{j,i}=
\frac{\gamma_{j,i}(l+q+2)}{q+2}.
\label{2.50}    
\end{equation}

With (\ref{2.50}) equality (\ref{2.49})
 takes the form
$$
{\cal A}(g_t)=
g_t+
\frac{(-1)^n}{q+2}
\sum_{r\in Q\backslash\{i\}}
\gamma_{j,r}
\sum_{p\in P\backslash\{j\}}
g_{r,p}\beta_{r,p}\eta_{ij}+
$$

\begin{equation}
{\hat g}_i(-1)^n
\gamma_{j,i}
\frac{l+q+2}{q+2}
\eta_{ij}.
\label{2.51}    
\end{equation}

Note that (\ref{2.34}) is correct as above. 
Then comparison of (\ref{2.45}) and  (\ref{2.51}) 
shows that to complete the proof of the lemma
it  suffices to check the equality 

$$
\frac{\hat\alpha_i}{\alpha_t}=
\frac{(-1)^n(q+1)}{q+2}.
$$

But the last one is correct for the same reason as 
  (\ref{2.34}) is.

Lemma 2.2 is completely proved.

Suppose $A_{h^{\prime}}$ and $A_g$ are as above.

{\bf Lemma 2.3.} {\em  If the matrix $A$ of
${\cal A}\in GL\; V$ has the form (\ref{2.7}),
then 
$$
A_{h^{\prime}}=A_g.
$$}

For any matrix $A$ of the form (\ref{2.7})
we have $\Ber A=1$. Hence, according to
Proposition 2.1, Lemma 2.3 is the next
part of the proof of Theorem 2.1 
(see (\ref{2.21})).

{\sc Proof of Lemma 2.3.} By assumption, there exist
$i\in \{1,\ldots,m\}$, $j\in \{1,\ldots,n\}$
such that

$$
{\cal A}(\varepsilon_j)=
\varepsilon_j+e_i\xi_{ij},
$$
where $\xi_{ij}\in G_1$, and for any
$w\in\Omega\backslash\{\varepsilon_j\}$ 
we have ${\cal A}(w)=w$.

Let $h_t$ be an element of $\Lambda_h$. By $R$ denote the
 subset of $\{1,\ldots,m\}$ such that $r\in R$ 
if and only if $\varepsilon_j$ belongs to the box
 $(r,j)$  of the $h_t$-tableau. Then we have

\begin{equation}
{\cal A}(h_t)=h_t+\sum_{r\in R}
(-1)^{n_{r,j}}{\tilde h}_r\xi_{ij},
\label{2.54}    
\end{equation}
where ${\tilde h}_r$ is obtained from $h_t$
using the replacement of $\varepsilon_j$ in the
 box $(r,j)$ by $e_i$;
as above,
$n_{r,j}$  is the number of odd 
base vectors that are placed from the box $(r,j+1)$ 
to the box $(m,n+1)$ when we move from left to right
and from the top down by the rows of the $h_t$-tableau.

Suppose $i\in R$, that is, the $h_t$-box $(i,j)$ is 
filled with
$\varepsilon_j$. Since $h_t$ is canonical,
we see that the $j$-th $h_t$-column does not 
contain $e_i$.

Evidently,  the tensor ${\tilde h}_i$ is canonical.
We claim that it is the only canonical tensor
among ${\tilde h}_r$, where $r\in R$. Indeed,
consider the $j$-th ${\tilde h}_r$-column,
where $r\neq i$. Then $e_i$ is placed not in its
``native''   $i$-th row but in the $r$-th row.

To reduce these tensors to a canonical form
we use their skew-symmetry by column elements
(see (\ref{1.3})). To be precise for a given $r\in R$
we transpose the vectors $e_i$ and $\varepsilon_j$
placed in the ${\tilde h}_r$-boxes $(i,j)$ and
$(r,j)$ respectively.

 As a result of this transposition
the additional factor $(-1)^{m_r+1}$  appears,
where $m_r$ is the number of odd base vectors
placed strictly between the $h_t$-boxes $(r,j)$ and
$(i,j)$ trough the movement by the rows from left to 
right and from top to bottom.

Since
$$
(-1)^{m_r+1}=(-1)^{n_{i,j}-n_{r,j}}
$$
we see that all terms of the sum in the right-hand 
side of (\ref{2.54}) are equal to
 $(-1)^{n_{i,j}}{\tilde h}_i\xi_{ij}$ and 
the number of these terms is $n_R=|R|$.
Then equality (\ref{2.54}) takes the form

$$
{\cal A}(h_t)=h_t+(-1)^{n_{i,j}}
n_R{\tilde h}_i\xi_{ij}.
$$

Going to the ``primed'' tensors we obtain

\begin{equation}
{\cal A}(h^{\prime}_t)=h^{\prime}_t+
(-1)^{n_{i,j}}n_R{{\tilde h}}^{\prime}_i
\frac{\alpha({\tilde h}_i)}{\alpha(h_t)}
\xi_{ij}.
\label{2.57}    
\end{equation}

If $i\notin R$ then the $h_t$-box $(i,j)$ is filled 
with $e_i$. As a result all tensors ${\tilde h}_r$,
where $r\in R$, are equal to zero because the element
$e_i$ enters twice in theirs $j$-th column but the
 tensors are skew-symmetric by elements of an 
arbitrary column. Thus for the case when $i\notin R$
we have
\begin{equation}
{\cal A}(h_t)=h_t.
\label{2.58}    
\end{equation}

Now let $g_t\in \Lambda_g$ be the tensor corresponding
 to  $h_t\in \Lambda_h$ considered above, that is, 
the tableaux of $h_t$ and $g_t$ have the same filling
within the small rectangle. By $\bar R$ denote the set
 of  integers such that $r\in \bar R$ if and only if
the $g_t$-box $(r,j)$ is filled with $\varepsilon_j$.
We claim that ${\bar R}=R\cup \{n+1 \}$. In fact, 
since the tensor $g_t$ is canonical, we see that
the $g_t$-box $(m+1,j)$ is filled with $\varepsilon_j$
and the tableaux of $h_t$ and $g_t$ are filled equally 
within the small rectangle.

Then we have
\begin{equation}
{\cal A}(g_t)=g_t+
\sum_{r\in \bar R}
(-1)^{{\bar n}_{r,j}}{\tilde g}_r\xi_{ij},
\label{2.59}    
\end{equation}
where ${\tilde g}_r$ is obtained from $g_t$ 
by the replacement of the element $\varepsilon_j$
disposed in the $g_t$-box $(r,j)$ to $e_i$; 
${\bar n}_{r,j}$ is the
number of odd base vectors placed in the $g_t$-tableau 
strictly after the box $(r,j)$ when we move
by the rows from left to right and from top to
 bottom.

Evidently, we have ${\bar n}_{r,j}=n_{r,j}+n$ for 
$r\leq m$ and ${\bar n}_{m+1,j}=n-j$.
Then separating the term for $r=m+1$ 
in right-hand side of 
(\ref{2.59})
 we obtain

\begin{equation}
{\cal A}(g_t)=g_t+
\sum_{r\in R}
(-1)^{n_{r,j}+n}{\tilde g}_r\xi_{ij}+
(-1)^{n-j}{\tilde g}_{m+1}\xi_{ij}.
\label{2.60}    
\end{equation}

Suppose that $i\in \bar R$, that is,
 the element $\varepsilon_j$ is placed
in the $g_t$-box $(i,j)$ and the element $e_j$
does not enter in the $j$-th $g_t$-column.

The tensor ${\tilde g}_i$ is canonical and the
 rest tensors ${\tilde g}_r$, where 
$r\in {\bar R}\backslash\{i\}$, are not canonical 
because the element $e_i$ 
is not placed in its ``native'' $i$-th row for 
these tensors.

To reduce a tensor ${\tilde g}_r$, where
$r\in {\bar R}\backslash\{i\}$, to a canonical form
we use  skew-symmetry
of the tensor by  elements of the column
(see (\ref{1.3})).

Then for any $r\in {\bar R}\backslash \{i\}$ we obtain 
the same tensor ${\tilde g}_i$ with the additional
 factor $(-1)^{{\bar m}_r+1}$, where ${\bar m}_r$ 
is the number of odd base vectors placed strictly
 between the boxes $(r,j)$ and $(i,j)$ (as above, 
we move from one box to another by the rows from
 left to right and from top to bottom). Since

$$
(-1)^{{\bar m}_r+1}=(-1)^{n_{r,j}-n_{i,j}}
$$
for $r\leq m$ and
$$
(-1)^{{\bar m}_{m+1}+1}=(-1)^{n_{i,j}+j},
$$ 
we see that all terms (distinct from $g_t$)
 in the right-hand side of
 (\ref{2.60}) are equal to
 $(-1)^{n_{i,j}+n}{\tilde g}_i\xi_{ij}$ and
the number of these terms is equal to 
$|{\bar R}|=n_R+1$.

Then equality (\ref{2.59}) takes the form
\begin{equation}
{\cal A}(g_t)=g_t+
(-1)^{n_{i,j}+n}(n_R+1)
{\tilde g}_i
\xi_{ij}.
\label{2.61}    
\end{equation}

Comparison of (\ref{2.57}) and (\ref{2.61}) 
shows that
to complete the proof of the lemma it  suffices
 to check the equality
\begin{equation}
\frac{\alpha({\tilde h}_i)}{\alpha(h_t)}=
\frac{(-1)^n(n_R+1)}{n_R}.
\label{2.62}    
\end{equation}
By definition (see (\ref{2.2})) we have
$$
\frac{\alpha({\tilde h}_i)}{\alpha(h_t)}=
\frac{(-1)^{n\rho({\tilde h}_i)}
\kappa({\tilde h}_i)\kappa(g_t)}
{\kappa({\tilde g}_i)(-1)^{n\rho(h_t)}\kappa(h_t)}.
$$
When passing from $h_t$ to ${\tilde h}_i$ the 
number of odd base vectors within the small
rectangle increases on $1$, that is,
$\rho(h_t)-\rho(h_i)=1$. 

The numbers of odd base vectors are the
same for all respective columns of  $h_t$ and
${\tilde h}_i$ ($g_t$ and ${\tilde g}_i$) except
the $j$-th one. Also, the numbers of odd base
vectors in the $j$-th column of $h_t$, ${\tilde h}_i$,
$g_t$, ${\tilde g}_i$ are equal to $n_R$,
$n_R-1$, $n_R+1$, $n_R$ (respectively).

Hence we have
$$
\frac{\kappa({\tilde h}_i)\kappa(g_t)}{
\kappa({\tilde g}_i)\kappa(h_t)}=
\frac{(n_R-1)!(n_R+1)!}{(n_R!)^2}=
\frac{n_R+1}{n_R}.
$$
and  equality (\ref{2.62})
is proved.

Suppose that $i\in {\bar R}$, that is, the box $(i,j)$ 
is filled with $e_i$. By the same argument as above
we get ${\cal A}(g_t)=g_t$.
Also from (\ref{2.58}) it follows that
${\cal A}(h^{\prime}_t)=h_t^{\prime}$.
Thus the case when $i\in {\bar R}$ is trivial.

This concludes the proof of  Lemma 2.3. 

A further step in the proof of Theorem 2.1
is the following lemma. 

{\bf Lemma 2.4.} {\em  If the matrix $A$ of
${\cal A}\in GL\; V$ has the form (\ref{2.10}),
then
$$
A_{h^{\prime}}=A_g.
$$ 
}

{\sc Proof of Lemma 2.4.} By assumption,
there exist unequal integers $s,t\in \{1,\ldots,n\}$
such that 
$$
{\cal A}(\varepsilon_s)=\varepsilon_s+
\varepsilon_ty,
$$
where $y\in G_0$, and for any
 $w\in \Omega\backslash\{\varepsilon_s\}$
we have
${\cal A}(w)=w$.

Let $h$ be an element of $\Lambda_h$.
By $F$ denote the set of integers such that
an integer $r$ belongs to $F$ if and only if the
 $h$-box $(r,s)$ is filled with $\varepsilon_s$.
 Suppose $F\neq \emptyset$.
Denote
$$
F_j=\{\nu|\;\nu\subseteq F,\;
|\nu|=j\},
$$
where $j=0,1,\ldots ,n_F$, $n_F=|F|$,
that is, $F_j$ is the set of all subsets
of $F$ with $j$ elements. Then we have 

\begin{equation}
{\cal A}(h)=h+\sum_{j=1}^{n_F}
\sum_{\nu\in F_j}{h_{j,\nu}}y^j,
\label{2.63}    
\end{equation}
where $h_{j,\nu}$ is obtained from $h$
using the replacement of $\varepsilon_s$
placed in the $h$-boxes $(i,s)$, where
$i\in \nu$, by $\varepsilon_t$.

For any $j>0$ an arbitrary element $h_{j,\nu}$
is not canonical, because $j$ elements of the form 
$\varepsilon_t$ are placed not in their ``native'' 
$t$-th column but in the $s$-th column. 
For an arbitrary  $h_{j,\nu}$ by 
$P(h_{j,\nu})$
denote the set of integers such that 
$p\in P(h_{j,\nu})$ if and only if $e_p$ 
belongs
to the $t$-th $h_{j,\nu}$-column and $e_p$ 
does not belong to
the $s$-th $h_{j,\nu}$-column. 
In other words, 
$p\in P(h_{j,\nu})$ iff the $h_{j,\nu}$-box $(p,t)$
is filled with $e_p$ but the $h_{j,\nu}$-box
 $(p,s)$
is filled  either with $\varepsilon_s$ or with
$\varepsilon_t$. 
Notice that for all
$j>0$, $\nu\in F_j$ the set $P(h_{j,\nu})$ 
is the
 same, that is, the composition of 
$P(h_{j,\nu})$
does not depend on $j$ and $\nu$.
For this reason in what follows we write simply $P$
instead of $P(h_{j,\nu})$. 

Denote
$$
P_j=\{\mu|\;\mu\subseteq P,\;|\mu|=j
\},
$$
where $j=0,1,\ldots ,n_F$, that is,  $P_j$
is the set of all subsets of $P$ with $j$ elements.

We claim that 
\begin{equation}
h_{j,\nu}=
\frac{1}{C_{q+j}^j}
\sum_{\mu\in P_j}
\delta_{\mu}{\bar h}_{j,\mu},
\label{2.64}    
\end{equation}
where $q$ is the number of odd base vectors
(i.e., the number of $\varepsilon_t$)
in the $t$-th column of $h$ and
${\bar h}_{j,\mu}$ is obtained from $h_{j,\nu}$
as follows: for every $z\in \mu$ the elements 
$e_z$ and $\varepsilon_t$ placed in the 
$h_{j,\nu}$-boxes $(z,t)$ and $(z,s)$ (respectively)
change their places;
$\delta_{\mu}=(-1)^{m_{\mu}}$, 
$m_{\mu}=\sum_{z\in {\mu}}m_z$,
where $m_z$ is the number of odd base vectors placed
in the $z$-th ${\bar h}_{j,\mu}$-row   strictly
 between the boxes $(z,t)$ and $(z,s)$.
Note that we can assume that the $h_{j,\nu}$-box
 $(z,s)$
is filled with $\varepsilon_t$. Indeed, 
otherwise the box
$(z,s)$ is filled with $\varepsilon_s$ and we use
(\ref{1.3}) to transpose $\varepsilon_s$ and
$\varepsilon_t$ in the $s$-th column.
Since the transposed elements are odd,
we see that the additional sign does not appear. 

Now let us prove identity (\ref{2.64}). To express
$h_{j,\nu}$ in terms of canonical tensors
we apply the Jacobi identity step by step.
In any step an element $\varepsilon_t$ from the $s$-th
column is moved to the $t$-th column.
Also the application of  the Jacobi identity in
 the next step brings the additional factor 
$1/(r+1)$, where $r$ is the number of elements 
$\varepsilon_t$ in the $t$-th column in the previous
step. Thus at last the factor 
$$
\frac{1}{(q+1)(q+2)\cdots(q+j)}
$$
arises.
Further, there is one-to-one correspondence 
between the set of canonical tensors obtained in
 the last step and the set of all subsets 
of $P$ with $j$ elements.
Also, any canonical tensor enters $j!$ times in the 
final expression for $h_{j,\nu}$,
because an element $e_z$, $z\in P$,
can come to the $s$-th column in any step of the 
algorithm and the number of these steps is equal 
to $j$.

Notice that the right-hand side of (\ref{2.64})
depends on $j$ and does not depend on $\nu$.
In other words, all tensors $h_{j,\nu}$,
where $\nu\in F_j$,
have the same representation in terms of 
the canonical tensors.
Then with (\ref{2.64}) equality (\ref{2.63}) 
takes the form 
$$
{\cal A}(h)=h+
\sum_{j=1}^{n_F}
\frac{C_{n_F}^j}{C_{q+j}^j}
\sum_{\mu\in P_j}
\delta_{\mu}{\bar h}_{j,\mu}
y^j
$$
or in terms of the ``primed'' tensors we have

\begin{equation}
{\cal A}(h^{\prime })= h^{\prime}+
\frac{1}{\alpha(h)}
\sum_{j=1}^{n_F}
\frac{C_{n_F}^j\alpha(h_{j,\mu})}{C_{q+j}^j}
\sum_{\mu\in P_j}
\delta_{\mu}{\bar h}^{\prime}_{j,\mu}
y^j,
\label{2.66}    
\end{equation}
where  $\displaystyle 
h^{\prime }=\frac{h}{\alpha(h)}$,
$\displaystyle {\bar h}^{\prime}_{j,\mu}=
\frac{h_{j,\mu}}{\alpha(h_{j,\mu})}$.

Now let $g\in \Lambda_g$ be the tensor 
corresponding to $h$ and ${\cal A}\in GL\; V$
as above. By $\bar F$ denote the set of integers 
such that $r\in \bar F$ if and only if the $g$-box
$(r,s)$ is filled with $\varepsilon_s$. Any column
 of $g$ as compared with $h$ contains 
the additional box filled with a proper odd base
 element (to be precise the $g$-box $(m+1,s)$
is filled with $\varepsilon_s$, where
 $s=1,2,\ldots,n$). Hence we have

$$
{\bar F}=F\cup \{m+1\}
$$
and $n_{\bar F}=|{\bar F}|=n_F+1$. Define

$$
{\bar F}_j=\{\nu|\;\nu\subseteq {\bar F},\; |\nu|=j\},
$$
where $j=0,1,\ldots,n_{\bar F}$,
that is, ${\bar F}_j$ is the set of all
 subsets of $\bar F$ with $j$ elements.
Then we have
\begin{equation}
{\cal A}(g)=g+
\sum_{j=1}^{n_F+1}
\sum_{\nu\in {\bar F}_j}
g_{j,\nu}y^j,
\label{2.67}    
\end{equation}
where $g_{j,\nu}$ are obtained from $g$
just as $h_{j,\nu}$ are obtained from $h$
(see above).

Since the tableaux of $h$ and $g$ have the same
 filling within the small rectangle and the
$(m+1)$-th $g$-row is filled with odd base vectors,
we see that $P(h_{j,\nu})=P(g_{j,\nu})$.
Therefore we write $P$ as above instead of 
$P(g_{j,\nu})$.

Let $P_j$ be as above. We see that $\delta_\mu$,
where $\mu\in P_{j}$, does not depend on either
$h_{j,\nu}$ or $g_{j,\nu}$ is considered.

We have 

\begin{equation}
g_{j,\nu}=\frac{1}{C_{q+j+1}^j}
\sum_{\mu\in P_j}
\delta_{\mu} {\bar g}_{j,\mu},
\label{2.68}    
\end{equation}
where 
${\bar g}_{j,\mu}$ are obtained from $g_{j,\nu}$
just as ${\bar h}_{j,\mu}$ are obtained from 
$h_{j,\nu}$.
Indeed, by analogy with (\ref{2.64}) equality 
(\ref{2.68})
can be obtained and the only difference between 
these two cases is that $g_{j,\nu}$
has not $q$ but $(q+1)$ elements $\varepsilon_t$
in the $t$-th column. That is why the coefficient in
(\ref{2.68}) is equal not to $1/C_{q+j}^j$ 
(as in (\ref{2.64})) 
but to $1/C_{q+j+1}^j$.

With (\ref{2.68}) equality (\ref{2.67}) takes
 the form
\begin{equation}
{\cal A}(g)=g+
\sum_{j=1}^{n_F}
\frac{C_{n_F+1}^j}{C_{q+j+1}^j}
\sum_{\mu\in P_j}
\delta_{\mu}
{\bar g}_{j,\mu}
y^j+
g_{\bar F}y^{n_F+1},
\label{2.69}    
\end{equation}
where $g_{\bar F}$ is obtained from $g$ using the 
replacement 
of {\it all}  $\varepsilon_s$
in the $s$-th $g$-column  by  $\varepsilon_t$.
Let us show that 
\begin{equation}
g_{\bar F}=0.
\label{2.70}    
\end{equation}
Indeed, applying the Jacobi identity we represent
 $g_{\bar F}$  as a linear combination of 
tensors such that theirs $s$-th column is filled 
with even base vectors only. But the height of the
 $s$-th column is greater by 1 than the number
$m$ of even base vectors. Since any tensor of the
 form (\ref{1.2}) is skew-symmetric 
 by elements of any column we get (\ref{2.70}).

Comparison of (\ref{2.66}) and (\ref{2.69})
shows that  to complete the proof it suffices
 to check the equality
$$
\frac{\alpha({\bar h}_{j,\mu})}{\alpha(h)}=
\frac{C_{q+j}^jC_{n_F+1}^j}{C_{n_F}^jC_{q+j+1}^j}
$$
or, what is the same -- 
\begin{equation}
\frac{\alpha({\bar h}_{j,\mu})}{\alpha(h)}=
\frac{(n_F+1)(q+1)}{(n_F-j+1)(q+j+1)}.
\label{2.71}    
\end{equation}

We have:

\begin{center}
\begin{tabular}{|c|c|c|}
\hline
 & \multicolumn{2}
{|c|}{the numbers of odd base vectors}\\
\cline{2-3}
 &{$s$-th column}&{$t$-th column}\\
\hline
{$h$}&{$n_F$}&{$q$}\\
\hline
{$g$}&{$n_F+1$}&{$q+1$}\\
\hline
{${\bar h}_{j,\nu}$}&{$n_F-j$}&{$q+j$}\\
\hline
{${\bar g}_{j,\nu}$}&{$n_F-j+1$}&{$q+j+1$}\\
\hline
\end{tabular}
\end{center}

Also note that the number of odd base vectors
within the small rectangle is the same for the
 tableaux of $h$ and ${\bar h}_{j,\mu}$,
that is, we have $\rho({\bar h}_{j,\mu})=\rho(h)$.
Hence we obtain
$$
\frac{\alpha({\bar h}_{j,\mu})}{\alpha(h)}=
\frac{\kappa(g)\kappa({\bar h}_{j,\mu})}
{\kappa(h)\kappa({\bar g}_{j,\mu})}=
$$
$$
=\frac{(n_F+1)!(q+1)!(n_F-j)!(q+j)!}
{n_F!q!(n_F-j+1)!(q+j+1)!}
$$
and we see that (\ref{2.71}) is proved.

This concludes the proof of Lemma 2.4.

Actually the following lemma is the last step
in the proof of Theorem 2.1.

{\bf Lemma 2.5.} {\em  If the matrix $A$ of
${\cal A}\in GL\; V$ has the form (\ref{2.9}),
then
$$
A_{h^{\prime}}=A_g.
$$ 
}

{\sc Proof.} By assumption, there exist integers
$s,t\in \{1,2,\ldots,m\}$, $s\neq t$ and $x\in G_0$
such that 
$$
{\cal A}(e_s)=e_s+e_tx,
$$
and for any $w\in \Omega\backslash \{e_s\}$
we have
${\cal A}(w)=w$.

Let $h$ be an element  of the set
 $\Lambda_h$. By $H$ denote the set of integers
 such that $r\in H$ if and only if the $h$-boxes
$(s,r)$ and $(t,r)$ are filled with $e_s$ and 
$\varepsilon_r$ (respectively). Evidently we have
$(n+1)\notin H$. By $H_j$ denote the set of all
 subsets of $H$ with $j$ elements, where 
$j=1,2,\ldots,k_H$, $k_H=|H|$. 

We have 
\begin{equation}
{\cal A}(h)=h+
\sum_{j=1}^{k_H}
\sum_{\mu\in H_j}
{h}_{\mu}x^j,
\label{2.72}    
\end{equation}
where ${h}_{\mu}$ is obtained from 
$h$ when all elements $e_s$ placed in the 
$h$-boxes $(s,r)$, where $r\in \mu$, are  
replaced by $e_t$.

An arbitrary tensor $h_{\mu}$ is not canonical
in so far as the $h_\mu$-boxes $(s,r)$,
where $r\in \mu$, are filled with $e_t$ but 
$t\neq s$.

To reduce a tensor ${h}_{\mu}$,
where $\mu\in H_j$, to a canonical form we use
the skew-symmetry of this tensor by the elements
 of the columns. To be precise,
for every $r\in \mu$ we transpose the vectors
$e_t$ and $\varepsilon_r$ placed in the $h_\mu$-boxes
$(s,r)$ and $(t,r)$ respectively.
The sign that appears as a result of 
this permutation  we denote 
by $\sigma_\mu$.

The canonical tensor obtained from ${h}_{\mu}$
is denoted by ${\bar h}_{\mu}$. 
Then equality (\ref{2.72})
takes the form

$$
{\cal A}(h)=h+
\sum_{j=1}^{k_H}
\sum_{\mu\in H_j}
\sigma_{\mu}{\bar h}_{\mu}x^j,
$$
and in terms of the ``primed'' tensors we have

\begin{equation}
{\cal A}({h^{\prime}})={h^{\prime}}+
\frac{1}{\alpha(h)}
\sum_{j=1}^{k_H}
\sum_{\mu\in H_j}
\sigma_{\mu}{{\bar h}^{\prime}}_{\mu}
\alpha({\bar h}_{\mu})
x^j,
\label{2.74}    
\end{equation}
where $h^{\prime}=h/\alpha(h)$, 
${\bar h}^{\prime}_{\mu}=
{\bar h}_{\mu}/\alpha(h_{\mu})$.

Note that for an arbitrary column the number of odd
base vectors is the same for the tableaux
of the tensors $h$, $h_\mu$ and ${\bar h}_\mu$.
Hence we have 
$\alpha({\bar h}_{\mu})=\alpha(h)$, 
where $\mu\in H_j$, $j=1,2,\ldots,k_H$.
Thus (\ref{2.74}) becomes

\begin{equation}
{\cal A}(h^{\prime})=h^{\prime}+
\sum_{j=1}^{k_H}
\sum_{\mu\in H_j}
\sigma_{\mu}{\bar h}^{\prime}_{\mu}x^j.
\label{2.75}    
\end{equation}

Now suppose $g$ is the tensor corresponding to 
$h$, that is, $g$ and $h$ have the same 
filling within the small rectangle.
Evidently the set $H$ is the same for the tensors
 $h$ and $g$. Then we get

$$
{\cal A}(g)=g+
\sum_{j=1}^{k_H}
\sum_{\mu\in H_j}
{g}_{\mu}x^j,
$$
where $g_\mu$ is obtained from $g$ much as
$h_\mu$ is obtained from $h$.
The tensors ${g}_{\mu}$ are reduced to a
 canonical form by the same process as 
${h}_{\mu}$ are. Then we obtain

\begin{equation}
{\cal A}(g)=g+
\sum_{j=1}^{k_H}
\sum_{\mu\in H_j}
{\hat \sigma}_{\mu}{\bar g}_{\mu}x^j,
\label{2.77}    
\end{equation}
where ${\bar g}_{\mu}$ are canonical tensors and we 
claim that ${\hat \sigma}_{\mu}=\sigma_\mu$. 
Indeed, $g_\mu$-boxes touched by the permutation
that reduces 
$g_\mu$ to a canonical form are 
placed strictly within the small
rectangle.

Comparison of (\ref{2.75}) and (\ref{2.77}) shows that
Lemma 2.5 is proved.

The case when the matrix $A$ of $\cal A$ has the form
(\ref{2.8}) is trivial.

Thus we see that Theorem 2.1 is completely proved.

Let $b_*$ be a formal element such that for 
any ${\cal A}\in GL\,V$ we have

$$
{\cal A}(b_*)=(\Ber A)^{-1}\cdot b_*,
$$
where $A$ is the matrix of $\cal A$.

{\bf Theorem 2.2.} {\em  The mapping 
$$
\varphi_*:b_*\cdot h_i^{\prime}\mapsto g_i
$$
determines the isomorphism of $GL\, V$-modules
$b_*\cdot W_{\lambda_h}$ and $W_{\lambda_g}$.  
}

{\sc Proof.} Follows immediately from Theorem 2.1.

\section{\rm An explicit construction of the 
one-dimensional representations of the general
 linear supergroup}
\setcounter{equation}{0}
\indent

By definition, put $\vartheta_i=e_i$ and 
$\vartheta_{j}
=\varepsilon_{j-m}$,
 where  $i=1,2,\ldots,m$,
$j=m+1,m+2,\ldots,m+n$. The last notation gives
 the through enumeration  of the set 
$\Omega$. 

By $V^*$ denote the $G$-module dual to $V$.

Let $l$ be a positive integer.
 
Suppose $u_i^*$, $w_j$ are arbitrary homogeneous
 elements of $V^*$ and $V$ (respectively),
where $i,j=1,2,\ldots,l$.
By definition, for the 
decomposable tensors
$u_1^*\cdots u_l^*\in T_l(V^*)$,
$w_1\cdots w_l\in T_l(V)$ put
\begin{equation}
(u_1^*\cdots u_l^*,w_1\cdots w_l)=
(-1)^\chi(u_1^*,w_1)\cdots(u_l^*,w_l),
\label{3.1}    
\end{equation}
where $\chi$ is the number of pairs of odd elements
that change their order when one pass from the 
sequence 
$u_1^*,\ldots,u_l^*,w_1,\ldots,$ $w_l$ to the sequence
$u_1^*,w_1,\ldots,$ $u_l^*,w_l$.
Equality (\ref{3.1}) determines the inclusion 
$$
T_l(V^*)\rightarrow(T_l(V))^*.
$$

By definition, for arbitrary 
$u_1^*,\ldots,u_l^*\in V^*$,  $\sigma\in S_l$ put
\begin{equation}
(u_1^*\cdots u_l^*) \sigma
=u_{\sigma(1)}^*\cdots u_{\sigma(l)}^*.
\label{3.2}    
\end{equation}
Thus $T_l(V^*)$
is the right $S_l$-module. With (\ref{1.1a}),
(\ref{3.1}), (\ref{3.2}) we get
\begin{equation}
((u_1^*\cdots u_l^*)\sigma,w_1\cdots w_l)=
(u_1^*\cdots u_l^*,\sigma(w_1\cdots w_l)),
\label{3.3}    
\end{equation}
where as above $u_i^*\in V^*$, $w_j\in V$.

By $\vartheta_i^*$ denote the elements of $V^*$
dual to $\vartheta_j$, that is,
$$
(\vartheta_i^*,\vartheta_j)=\delta_{ij},
$$
where $i,j=1,2,\ldots,m+n$,
$\delta_{ij}$ is the Kronecker delta. 
Clearly, $e_i^*=\vartheta_i^*$,
$\varepsilon_j^*=\vartheta_{j+m}^*$,
where $i=1,\ldots,m$, $j=1,\ldots,n$. 
By analogy 
with the notation used above, 
$\Omega_0^*=\{e_1^*,\ldots,e_m^*\}$,
$\Omega_1^*=\{\varepsilon_1^*,\ldots,\varepsilon_n^*\}$,
$\Omega^*=\Omega_0^*\cup\Omega_1^*$.

We say that $f\in T_l(V^*)$ is 
an $\Omega^*$-{\it tensor} 
if there exist 
$u_1^*,\ldots,u_l^*\in \Omega^*$  
and a Young tableau $T$
such that

$$
f=(u_1^*\cdots u_l^*) e_T.
$$

Let $\Lambda_g$ be as above (see \S 1), $g_i$  
elements of $\Lambda_g$ 
in some ordering of the last one, where 
$i=1,2,\ldots,k_\Lambda$; $k_\Lambda=2^{mn}$.
 For an arbitrary $i\in\{1,\ldots,k_\Lambda\}$
by $g_i^*$ denote the $\Omega^*$-tensor such
that $g_i^*$-tableau is obtained from $g_i$-tableau
 by the formal change
$\vartheta_j\mapsto \vartheta_j^*$, where 
$j=1,2,\ldots,(m+n)$.

By $k_i$ denote the cardinality of the 
automorphism group of the $g_i$-tableau and by $q_i$ 
denote the number of odd base elements in the 
$g_i$-tableau.

{\bf Lemma 3.1.} {\em  For arbitrary 
$i,j\in \{1,2,\ldots,k_\Lambda\}$
the following equality holds
\begin{equation}
(g_i^*,g_j)=\delta_{ij}\zeta_i,
\label{3.4}    
\end{equation}
where  
$$
\zeta_i={k_i}
\mu_{\lambda_g}(-1)^{(q_i^2-q_i)/2},
$$

\begin{equation}
\mu_{\lambda_g}=
\frac{\displaystyle\prod_{t=1}^{m+1}(m+n+1-t)!}
{\displaystyle\prod_{t=1}^{m+1}(m+1-t)!}.
\label{3.4.1}    
\end{equation}

} 

{\sc Proof.} By definition,
there exist $u_1^*,\ldots,u_l^*\in \Omega^*$,
$w_1,\ldots,w_l\in \Omega$ such that
$$
g_i^*=(u_1^*\cdots u_l^*)e_T,\;\;\;
g_j=e_T(w_1\cdots w_l),
$$
where $T$ is the Young tableau with the diagram
$\lambda_g=(\underbrace{n,\ldots,n}_{m+1})$,
$l=(m+1)n$.

With (\ref{3.3}) we get
$$
(g_i^*,g_j)=
((u_1^*\cdots u_l^*)e_T,e_T(w_1\cdots w_l))=
(u_1^*\cdots u_l^*,e_T^2(w_1\cdots w_l)).
$$

For an arbitrary Young diagram $\lambda$ 
the following identity holds
$$
e_{T_\lambda}^2=\mu_\lambda e_{T_\lambda},
$$
where $T_\lambda$ is a Young tableau with the 
diagram $\lambda$,
$\mu_\lambda$ is a non-zero integer. In particular, 
$\mu_{\lambda_g}$ is given by (\ref{3.4.1}) 
(see \cite{Weyl}). By assumption,
$u_1^*,\ldots,u_l^*$, $w_1,\ldots,w_l$ fill the diagram $\lambda_g$
in a canonical way. Then to conclude the proof 
it suffices to use Lemma 1.1 from \cite{Det}.

Thus Lemma 3.1 is proved.

An arbitrary linear transformation
 ${\cal A}\in GL\,V$ determines, as it usually 
does, the corresponding linear transformation 
of $V^*$  by the equalities
\begin{equation}
({\cal A}(\vartheta_i^*),{\cal A}(\vartheta_j))=
(\vartheta_i^*,\vartheta_j)=\delta_{ij}
\label{3.5}    
\end{equation}
(this transformation of $V^*$ we still  
denote by $\cal A$).

Suppose $A$ is the matrix of $\cal A$, that is, 
we have
\begin{equation}
({\cal A}(\vartheta_1),\ldots,{\cal A}(\vartheta_{m+n}))=
(\vartheta_1,\ldots,\vartheta_{m+n})\cdot A.
\label{3.6}    
\end{equation}

Then with (\ref{3.5}) we obtain
\begin{equation}
\left( 
\begin {array}{c}
{\cal A}(\vartheta_1^*) 
\\
{\cal A}(\vartheta_2^*)
\\
\cdots
\\
{\cal A}(\vartheta_{m+n}^*)
\end{array} 
\right)=
A^{-1}\cdot
\left( 
\begin {array}{c}
\vartheta_1^* 
\\
\vartheta_2^*
\\
\cdots
\\
\vartheta_{m+n}^*
\end{array} 
\right).
\label{3.7}    
\end{equation}

Denote
$$
g_i^{*\prime}=
\frac{1}{(g_i^*,g_i)}
g_i^*.
$$

From (\ref{3.4}) it follows that
$$
(g_i^{*\prime},g_j)=\delta_{ij}.
$$

Recall that $\cal A$ is even. Then it follows from
(\ref{3.1}) and (\ref{3.5}) that 

$$
({\cal A}(g_i^{*\prime}),{\cal A}(g_j))=
(g_i^{*\prime},g_j)=\delta_{ij}.
$$
In other words, the action of $\cal A$ on
$ \langle g_i \rangle$ and 
$\langle g_i^{*\prime}\rangle$ is the particular case
of changing base in a space and the dual one
(see (\ref{3.6}), (\ref{3.7})).

Consequently we have
$$
\sum_{i=1}^{k_\Lambda}
{\cal A}(g_i){\cal A}(g_i^{*\prime})=
({\cal A}(g_1),\ldots,{\cal A}(g_{k_\Lambda}))
\left( 
\begin {array}{c}
{\cal A}(g_1^{*\prime})
\\
\cdots
\\
{\cal A}(g_{k_\Lambda}^{*\prime})
\end{array} 
\right)=
$$
 
$$
=(g_1,\ldots,g_{k_\Lambda})
{\tilde A}{\tilde A}^{-1}
\left( 
\begin {array}{c}
g_1^{*\prime}
\\
\cdots
\\
g_{k_\Lambda}^{*\prime}
\end{array} 
\right)=
(g_1,\ldots,g_{k_\Lambda})
\left( 
\begin {array}{c}
g_1^{*\prime}
\\
\cdots
\\
g_{k_\Lambda}^{*\prime}
\end{array} 
\right)=
\sum_{i=1}^{k_\Lambda}g_ig_i^{*\prime}.
$$

This completes the proof of the following result:

{\bf Proposition 3.1.}  {\em  The action 
of an arbitrary
${\cal A}\in GL\, V$ on the tensor 
$$
\sum_{i=1}^{k_\Lambda}
g_ig_i^{*\prime}
$$ 
is identical.}

{\bf Theorem 3.1.} 
{\em The tensor 
$$
{\tilde b}=\sum_{i=1}^{k_\Lambda}
h^{\prime}_ig_i^{*\prime}
$$
generates the one-dimensional $GL\, V$-module
such that
$$
{\cal A}({\tilde b})=\Ber A\cdot {\tilde b},
$$
 where $A$ is a matrix of  
${\cal A}\in GL\, V$.}

{\sc Proof.} Follows immediately from Theorem 2.1 and
Proposition 3.1.

Consider some examples.

{\bf Example 3.1.} Suppose $n=0$. In this case one
 can say that the small rectangle degenerates into
the segment of height $m$. Then ${\tilde b}
=\gamma h_1$, where $h_1$ is the tensor skew-symmetric
in $e_1,\ldots,e_m$; $\gamma$ is an 
invertible element of $G_0$. 
Also, for $n=0$ we have $\Ber A=\det A$.
Thus we arrive at the basic classic result.

{\bf Example 3.2.} Let the dimension of $V$
be $1|1$; $V_0=\langle e\rangle$, 
$V_1=\langle \varepsilon\rangle$.
By ${\cal T}_{h_i}$  (${\cal T}_{g_i}$)
denote the $h_i$-tableau ($g_i$-tableau, 
respectively). Then we have

\begin{center}
\begin{picture}(150,16)
\put(32,0){\line(1,0){32}}
\put(32,16){\line(1,0){32}}
\put(32,0){\line(0,1){16}}
\put(48,0){\line(0,1){16}}
\put(64,0){\line(0,1){16}}
\put(1,6){${\cal T}_{h_1}=$}
\put(37,6){$e$}
\put(53,6){$e$}
\put(66,8){,}

\put(112,0){\line(1,0){32}}
\put(112,16){\line(1,0){32}}
\put(112,0){\line(0,1){16}}
\put(128,0){\line(0,1){16}}
\put(144,0){\line(0,1){16}}
\put(117,6){$\varepsilon$}
\put(133,6){$e$}
\put(81,6){${\cal T}_{h_2}=$}
\put(146,8){,}
\end{picture}
\end{center}

\vspace{-0.5cm}

$$
h_1=({\bf e}+(12))ee=2ee,
$$
$$
h_2=({\bf e}+(12))\varepsilon e=
\varepsilon e+e\varepsilon,
$$
where $\bf e$ is the unit of $S_2$. With
(\ref{2.2}) and (\ref{2.3}) we get

\begin{equation}
h^{\prime}_1=2ee,\;\;\;
h^{\prime}_2=-2(\varepsilon e+e\varepsilon).
\label{3.14}    
\end{equation} 

Further, we have

\begin{center}
\begin{picture}(150,32)
\put(32,0){\line(1,0){16}}
\put(32,16){\line(1,0){16}}
\put(32,32){\line(1,0){16}}
\put(32,0){\line(0,1){32}}
\put(48,0){\line(0,1){32}}
\put(1,14){${\cal T}_{g_1}=$}
\put(37,6){$\varepsilon$}
\put(37,22){$e$}
\put(50,16){,}

\put(112,0){\line(1,0){16}}
\put(112,16){\line(1,0){16}}
\put(112,32){\line(1,0){16}}
\put(112,0){\line(0,1){32}}
\put(128,0){\line(0,1){32}}
\put(117,6){$\varepsilon$}
\put(117,22){$\varepsilon$}
\put(81,14){${\cal T}_{g_2}=$}
\put(130,16){,}
\end{picture}
\end{center}
$q_1=k_1=1$,  $q_2=k_2=2$, $\mu_{\lambda_g}=2$.
 Hence, $\zeta_1=2$, $\zeta_2=-4$.
Therefore,

\begin{equation}
g_1^{*\prime}=
\frac{1}{2}
(e^*\varepsilon^*({\bf e}-(12)))=
\frac{1}{2}
(e^*\varepsilon^*-\varepsilon^*e^*),
\label{3.15}    
\end{equation}

\begin{equation}
g_2^{*\prime}=
-\frac{1}{4}
(\varepsilon^*\varepsilon^*({\bf e}-(12)))=
-\frac{1}{2}
\varepsilon^*\varepsilon^*.
\label{3.16}    
\end{equation}

Now with (\ref{3.14}), (\ref{3.15}), (\ref{3.16})
we obtain
$$
{\tilde b}=
\sum_{i=1}^{2}h^{\prime}_ig^{*\prime}_i=
eee^*\varepsilon^*-
ee\varepsilon^*e^*+
e\varepsilon\varepsilon^*\varepsilon^*+
\varepsilon e\varepsilon^*\varepsilon^*.
$$

\vspace{0.5cm}

Let $\Lambda_h=\{h_i\}$ be as above (see \S 1).
For an arbitrary $i\in \{1,2,\ldots,k_\Lambda\}$
by $h_i^*$ denote the $\Omega$-tensor such that the
$h_i^*$-tableau is obtained from the $h_i$-tableau
by the formal change 
$\vartheta_j\mapsto\vartheta_j^*$, 
where $j=1,2,\ldots,(m+n)$.

By $l_i$ denote the cardinality of the 
automorphism group of the $h_i$-tableau and
by $p_i$ denote the number of odd elements
 of the $h_i$-tableau.

Then by analogy with Lemma 3.1 we obtain:

{\bf Lemma 3.2.}  {\em  For arbitrary 
$i,j\in \{1,2,\ldots,k_\Lambda\}$ the following 
equality holds
$$
(h_i^*,h_j)=\delta_{ij}\zeta_i^{\prime},
$$
where
$$
\zeta_i^{\prime}={l_i}
\mu_{\lambda_h}(-1)^{(p_i^2-p_i)/2},
$$
$$
\mu_{\lambda_h}=
\frac{\displaystyle\prod_{t=1}^{m}(m+n+1-t)!}
{\displaystyle\prod_{t=1}^{m}(m-t)!}.
$$
}

By definition, put 

$$
h_i^{*\prime}=\frac{\alpha(h_i)}
{(h_i^*,h_i)}
h_i^*.
$$
With (\ref{2.3}) we have

$$
(h_i^{*\prime},h_j^{\prime})=\delta_{ij}.
$$

By analogy with Proposition 3.1 and
 Theorem 3.1 we arrive at the following results:

{\bf Proposition 3.2.} 
{\em The action of an arbitrary ${\cal A}\in GL\,V$
on the tensor 
$$
\sum_{i=1}^{k_\Lambda}
h_i^{\prime}h_i^{*\prime}
$$
is identical.}

{\bf Theorem 3.2.} 
{\em The tensor 
$$
{\tilde b}_*=\sum_{i=1}^{k_\Lambda}
g_ih_i^{*\prime}
$$
generates the one-dimensional $GL\, V$-module
such that
$$
{\cal A}({\tilde b}_*)=(\Ber A)^{-1}\cdot 
{\tilde b}_*,
$$
 where $A$ is a matrix of  
${\cal A}\in GL\, V$.}

{\bf Example 3.3.}  Let the dimension of $V$
be $1|1$; $V_0=\langle e\rangle$, 
$V_1=\langle \varepsilon\rangle$. Then we have

$$
h_1^*=e^*e^*({\bf e}+(12))=2e^*e^*,
$$

$$
h_2^*=\varepsilon^*e^*({\bf e}+(12))=
\varepsilon^*e^*+e^*\varepsilon^*
$$
(see the $h_i$-tableaux ${\cal T}_{h_i}$ 
in Example 3.2). Also, $l_1=2$, $p_1=0$,
$l_2=p_2=1$, $\mu_{\lambda_h}=2$.
Hence, $\zeta_1^{\prime}=4$,
$\zeta_2^{\prime}=2$.

With (\ref{2.3}) we get $\alpha(h_1)=1$,
${\displaystyle\alpha(h_2)=-\frac{1}{2}}$.
Thus we obtain

$$
h_1^{*\prime}=\frac{1}{2}e^*e^*,\;\;
h_2^{*\prime}=-\frac{1}{4}(\varepsilon^*e^*+
e^*\varepsilon^*).
$$

Further,
$$
g_1=({\bf e}-(12))e\varepsilon=
e\varepsilon-\varepsilon e,
$$

$$
g_2=({\bf e}-(12))\varepsilon\varepsilon=
2\varepsilon\varepsilon
$$
(see the $g_i$-tableaux ${\cal T}_{g_i}$ 
in Example 3.2).
At last, we obtain 

$$
{\tilde b}_*=\sum_{i=1}^{2}g_ih_i^{*\prime}
=
\frac{1}{2}
(e\varepsilon e^*e^*-
\varepsilon ee^*e^*-
\varepsilon\varepsilon\varepsilon^* e^*-
\varepsilon\varepsilon e^*\varepsilon^*).
$$

\end{document}